\documentclass[12pt]{amsart}
\usepackage{amsmath, amssymb, amsthm, latexsym}
\input xypic
\newtheorem{theorem}{Theorem}%[section]

\newtheorem{lemma}[theorem]{Lemma}

\newtheorem{corollary}[theorem]{Corollary}
\newtheorem{definition}[theorem]{Definition}

\def\Proof{\medskip\noindent{\bf Proof: }}

\def\Z{\mathbb{Z}}

\def\P{\mathbb{P}}

\def\A{\mathbb{A}}

\def\C{\mathbb{C}}

\def\Q{\mathbb{Q}}

\def\C{\mathbb{C}}

\def\Pi{\mathbb{P}^{\infty}}

\def\qed{\hfill$\square$\medskip}

\def\Zpk{\mathbb{Z}/p^{k}}
\def\Zpk1{\mathbb{Z}/p^{k-1}}

\newcommand{\rref}[1]{(\ref{#1})}

\newcommand{\cform}[3]{\begin{array}{c}
{\scriptstyle #3}\\
#1\\
{\scriptstyle #2}\end{array}}

\newcommand{\fracd}[2]{\frac{\displaystyle #1}{\displaystyle #2}}

\newcommand{\beg}[2]{\begin{equation}\label{#1}#2\end{equation}}
\def\r{\rightarrow}
\def\mc{\mathcal{C}}

\def\sl2{\widetilde{SL_{2}(\Z)}}

\title[Tree Algebras]{Tree algebras: An algebraic axiomatization of intertwining vertex
operators}
\author{Igor Kriz and Yang Xiu}
\thanks{The first author was supported in part by NSA grant H98230-09-1-0045. 
The second author was supported by a Fellowship from the University of Michigan
}
\begin{document}

\maketitle

\begin{abstract}
We describe a completely algebraic axiom system for intertwining operators
of vertex algebra modules, using algebraic flat connections, thus
formulating the concept of a {\em tree algebra}. Using 
the Riemann-Hilbert correspondence, we further prove that a vertex tensor
category in the sense of Huang and Lepowsky gives rise to a tree algebra over
$\C$. We also show that the chiral WZW model of a simply connected simple
compact Lie group gives rise to a tree algebra over $\Q$.
\end{abstract}

\section{Introduction}

The aim of the present paper is to describe a set of purely algebraic 
axioms
designed to capture the structure of
genus $0$ (or ``tree level'') amplitudes in conformal field
theory. To explain this very briefly, in mathematical physics,
a chiral quantum conformal field theory consists of a system of
state spaces indexed by a set of {\em labels} $\Lambda$. A Riemann
surface with parametrized boundary components whose boundary
components are labelled by elements of
$\Lambda$ should specify, roughly, a 
linear operator from the tensor product of the state
spaces corresponding to the inbound boundary components
to the tensor product of the state spaces corresponding to
the outbound boundary components. In fact, the operator
is not specified uniquely. Rather, there should be
a finite vector space of such operators, and these vector
spaces should form a holomorphic bundle on the moduli space
of Riemann surfaces, called a {\em modular functor}.
Certain axioms should
be satisfied, in particular, when cutting a Riemann surface $\Sigma$
along an analytically parametrized simple curve, the space of
operators corresponding to $\Sigma$ should be isomorphic,
via a prescribed isomorphism, to the direct sum
of the vector spaces corresponding to all possible labels on
the cutting curve (the spaces should be $0$ for all but finitely
many choices of labels on the cutting curve). The isomorphisms
should be subject to certain canonical coherence diagrams. The operators
on state spaces corresponding to the uncut and
cut Riemann surfaces should be related by trace. There should also
be a special label called the {\em zero label}
such that the operator corresponding to annuli approaching
the unit circle ``converge'' to the identity. 

\vspace{3mm}

The main issue with such an approach is inherited from the
usual issues of quantum mechanics: in interesting examples,
the state spaces are infinite-dimensional, and to make the
approach mathematically rigorous, we must select which category
of vector spaces
we will be working in. Quantum mechanics suggests Hilbert spaces
(cf. \cite{scft}),
in which case we are firmly in the realm of analysis.
Involving analysis certainly seems necessary if we want
to make rigorous the entire structure outlined in the last
paragraph.

\vspace{3mm}
In the 80's, however, it has been noticed by Borcherds \cite{bor}
and Frenkel, Lepowsky and Meurman \cite{flm} that the structure
present on the state space corresponding to the $0$ label
and genus $0$ Riemann surfaces
can be axiomatized purely algebraically. In this approach,
the state space is not a Hilbert space, but simply a graded
vector space over a field of characteristic $0$. Operators
corresponding to round domains (the unit disk with a finite
set of disks removed, all boundary components parametrized
linearly) are encoded in a structure called {\em (graded)
vertex algebra}. If one wishes to look beyond round domains,
one does not consider arbitrary Riemann surfaces of genus $0$,
but rather infinitesimal variations of boundary parametrization,
the effect of which is
described by what is known as the {\em conformal element}
or {\em energy-momentum tensor}.

\vspace{3mm}

The goal of the present paper is to propose an extension
of the vertex algebra axiomatization to {\em intertwining operators},
i.e. operators on state spaces involving the entire set of labels
$\Lambda$ of a chiral conformal field theory,
purely in the language of algebra. In particular,
a test of the success of such an approach is that it
should be valid over any field of characteristic $0$.
(Note: in this paper, we only deal with round domains,
we do not consider conformal elements). This advances
the program of making conformal field theory 
a purely algebraic object, which is a well known desideratum. It is
related, for example, to the 
Grothendieck-Teichm\"{u}ller program
\cite{esquisse,ideal}, and, viewed from a different perspective, the
geometric
Langlands program (cf. \cite{beidr}).

\vspace{3mm}
There is, however, a major difficulty. It is well known
that in interesting examples (e.g. parafermions or
the chiral WZW models), the intertwining operators 
corresponding to non-zero
labels (even for round domains) fail to be algebraic.
Typically, rather, one gets hypergeometric functions (which
are transcendental) and their generalizations (cf. \cite{gh,hl1}).
Considering this, our task may seem impossible. Because of
this, axiomatic systems have been developed which combine
algebra and analysis, in particular the concept of 
{\em vertex tensor category} by Huang and Lepowsky \cite{hl}.

\vspace{3mm}
Yet, the same examples point toward a possible solution:
although the correlation functions for non-zero labels are
not algebraic, they satisfy {\em differential equations} which
are algebraic. An example is given by the Knizhnik-Zamolodchikov
equations in the case of the WZW model (cf. \cite{hl1}). In
fact, similar equations have been formulated by
Y.Z. Huang for a large range
of examples in \cite{h7}. 

\vspace{3mm}
What we do in the present paper is formulate purely algebraic
axioms the Huang-Knizhnik-Zamolodchikov equations should
satisfy in an algebraic model of chiral conformal field
theory in genus $0$. The main obstacle to doing so is the
sheer complexity of the structure involved. Even the algebraic
axiomatization of a vertex algebra \cite{bor,flm} (i.e. the case of the
$0$ label) involves the remarkably complicated Jacobi identity. 
To make things worse,
a straightforward generalization of the Jacobi identity to intertwining
operators is false. Although Huang \cite{hjac} has obtained a 
more complicated generalization
of the Jacobi identity for intertwining operators, devicing a workable
system of axioms from this approach seems daunting.
To get a handle on the structure in the present paper, we make
substantial use a simple interpretation
of the vertex algebra axioms obtained by Hortsch, Kriz and Pultr \cite{hkp}:
a graded vertex algebra is essentially the same thing as an algebra over
a certain co-operad (the ``correlation function co-operad''), which can
be easily described explicitly. 

\vspace{3mm}
The present paper is organized as follows:
We write down the general axioms
in Section \ref{sdef} below; we call the resulting
algebraic structure a {\em tree algebra}. In Section \ref{sreg},
we discuss a version of the Riemann-Hilbert correspondence (cf. \cite{d})
in the present setting, and prove that a vertex tensor category
in the sense of Huang and Lepowsky \cite{hl} can be realized 
as a tree algebra over $\C$ ``with regular singularities'',
which is a completely algebraic object. 
At the end of Section \ref{sreg}, we also give an example showing
that the chiral WZW model for a simply connected simple compact
Lie group can be realized as a tree algebra over $\Q$.

\section{The basic definitions}
\label{sdef}

Let $F$ be any field of characteristic $0$. We recall from
\cite{hkp} the ordered configuration space $C(n)=C(z_1,...,z_n)$ of $n$ 
distinct points $z_1,...,z_n$ in $\A^1$. Denote the 
coefficient ring of the affine variety $C(n)$ by
$\mc(n)$. It is advantageous to consider $\mc(n)=\mc(z_1,...,z_n)$
a $\Z$-graded ring with grading by total degree of
a homogeneous rational function. In \cite{hkp}, it is
shown that the system $(\mc(n))$ is a graded co-operad:
the co-multiplication 
\beg{egeneral+}{\begin{array}{l}\mc(z_{11},...,z_{1n_1},...,z_{k1},...,z_{kn_k})_\ell\\
\r\mc(z_1,...,z_k)_{\ell_0}\otimes\mc(t_{11},...,t_{1n_1})_{\ell_1}
\otimes ...\otimes \mc(t_{k1},...,t_{kn_k})_{\ell_k}\end{array}
}
with $\ell_0+...+\ell_k=\ell$ is given by setting
$$z_{ij}=t_{ij}+z_i$$
and expanding
\beg{egeneral++}{(t_{ij}+z_i-t_{i^\prime j^\prime}-z_{i^\prime})^{-1}, \; i\neq i^\prime
}
by rewriting \rref{egeneral++} as
$$(t_{ij}+(z_i-z_{i\prime}) -t_{i^\prime j^\prime})^{-1},$$
and expanding in increasing powers of $t_{ij}$ and $t_{i^\prime j^\prime}$. 
Moreover, it is shown in \cite{hkp}
that vertex algebras in a suitable sense can be characterized as graded algebras
over the graded co-operad $\mc$.

Define, for $\Z$-graded $F$-vector spaces $A$, $B$, a $\Z$-graded vector space
$$(A\widehat{\otimes}B)_n=\cform{colim}{k_0}{}\cform{\prod}{k\geq k_0}{} 
A_k\otimes B_{n-k}$$ 
If $A$, $B$ are (commutative) rings, so is (in a natural way) $A\widehat{\otimes} B$. 
One checks that 
the structure map \rref{egeneral+} takes the form
\beg{egen+1}{\mc(m_1+...+m_n)\r (\mc(m_1)\otimes ...\otimes \mc(m_n))
\widehat{\otimes} \mc(n),}
and 
is a homomorphism of graded rings. 

Note that $\mc(0)=F$, so selecting $i$ among $n$ coordinates gives us a 
``co-augmentation''
\beg{eaug}{\mc(n-i)\r \mc(n).
}

\vspace{3mm}
The purpose of this paper is to extend this definition to a fully algebraic
treatment of tree-level amplitudes of chiral conformal field theories. 
This makes it possible to define, at least on the level of chiral tree-level
amplitudes, ``rational conformal field theory'' as ``conformal field theory
over $\Q$''.
This is possible despite of the fact that these amplitudes are usually
transcendental: typical examples are hypergeometric functions
\cite{gh}. The reason that an algebraic treatment is possible is because
we can describe {\em algebraic flat connections with regular singularities}
whose solutions are the desired amplitudes, and characterize
the precise algebraic structure
these connections are required to obey. To this end, we must first
develop the theory of such connections: the reason this is non-trivial
is that we need a suitable treatment of the grading, which in particular
should capture a concept of ``total degree'' of the solution, which is
to be an element of the field $F$.

\vspace{3mm}
We consider the ring of K\"{a}hler differentials $\Omega^{1}_{\mc(n)/F}$,
and equip it with a $\Z$-grading so that the universal differentiation
$$d:\mc(n)\r \Omega^{1}_{\mc(n)/F}$$
is a graded homomorphism of $F$-modules (of degree $0$). Then $\Omega^{1}_{\mc(n)/F}$
is a free $\mc(n)$-module on the basis $dz_1,...,dz_n$.
Let $M$ be a projective graded $\mc(n)$-module. Then 
a {\em homogeneous connection} is a map of graded $F$-modules
$$\nabla:M\r M\otimes_{\mc(n)}\Omega^{1}_{\mc(n)/F}
$$
which satisfies, for $f\in \mc(n)$, $m\in M$,
$$\nabla(fm)=mdf+f\nabla(m).
$$
We say that the connection $\nabla$ is {\em flat} when it satisfies the 
Maurer-Cartan equation
\beg{emc*}{(1\otimes d)\nabla=-(\nabla\otimes 1)\nabla
}
where both sides are considered as maps into $M\otimes_{\mc(n)}\Omega^{2}_{\mc(n)/F}$.

\vspace{3mm}
Next, our aim is to define the {\em degree} of a flat homogeneous
connection. We first define the degree of the {\em difference}
of two flat homogeneous connections. In effect, such difference
$$E=\nabla_1-\nabla_2,$$
written in matrix form as
\beg{emc1}{E=E_1dz_1+...+E_ndz_n,\; E_i\in Hom_{\mc(n)}(M,M)
}
is said to {\em have a degree} when
\beg{emc2}{E_1z_1+...+E_nz_n = k(z_1,...,z_n)Id_{M}
}
for some function
$$k(z_1,...,z_n)\in \mc(n).$$
Clearly, this property is invariant under linear change of the variables $z_1,...,z_n$.

\vspace{3mm}
\begin{lemma}
\label{lmc1}
When a difference $E$ of two flat connections has a degree, then $deg(E)=k(z_1,...,z_n)$
is a constant function of $z_1,...,z_n$.
\end{lemma}

\Proof
In coordinates, \rref{emc*} reads
\beg{emc3}{\frac{\partial{E_i}}{\partial z_j}-\frac{\partial{E_j}}{\partial z_i}
=-\frac{1}{2}[E_i,E_j].
}
Now compute
\beg{emc4}{\begin{array}{l}
\fracd{\partial(E_1z_1+...+E_nz_n)}{\partial z_1}=E_1+
z_1\fracd{\partial E_1}{\partial z_1} +...+ z_n \fracd{\partial E_n}{\partial z_1}=\\[3ex]
=E_1+z_1\fracd{\partial E_1}{\partial z_1}+...+z_n\fracd{\partial E_1}{\partial z_n}-
\fracd{1}{2}[E_1,z_1E_1+...+z_nE_n].
\end{array}
}
The Lie bracket \rref{emc4} vanishes because of our assumption \rref{emc2}.
The first term vanishes by the following Lemma and the fact that
$E$ is a homogeneous connection. Thus, \rref{emc4} vanishes. An analogous
argument holds with $z_1$ replaced by any $z_i$, which proves the
statement of the Lemma.
\qed

\vspace{3mm}

\begin{lemma}
\label{lmc2}
A function $f$ of $n$ variables $z_1,...,z_n$ is homogeneous
of degree $k$ if 
and only if
\beg{ehomog}{z_1\frac{\partial f}{z_1}+...+z_n\frac{\partial f}{\partial z_n}= kf.}
\end{lemma}

\Proof
Suppose $f$ is homogeneous of degree $k$.
By the chain rule,
$$\begin{array}{l}
ka^{k-1}f(z_1,...,z_n)=\fracd{da^k}{da}f(z_1,...,z_n)=\fracd{df(az_1,...,az_n)}{da}=\\[3ex]
=\left.\fracd{\partial f(t_1,...,t_n)}{\partial t_1}\right|_{(az_1,...,az_n)}\cdot z_1
+...+
\left.\fracd{\partial f(t_1,...,t_n)}{\partial t_n}\right|_{(az_1,...,az_n)}\cdot z_n=
\\[4ex]
=a^{k-1}\left(z_1\fracd{\partial f}{\partial z_1}+...+z_n\fracd{\partial f}{\partial z_n}\right).
\end{array}
$$
Conversely, if \rref{ehomog} holds,
$$\frac{df(\lambda z_1,...,\lambda z_n)}{d\lambda}=\left.\sum z_i\frac{\partial f(t_1,...,t_n)}{
\partial t_i}\right|_{t=\lambda z_i}=$$
$$=\frac{kf(\lambda z_1,...,\lambda z_n)}{\lambda}.$$
Studying the solutions of the ODE $y^\prime(\lambda)=\frac{ky}{\lambda}$ gives
the result.
\qed

\vspace{3mm}

\begin{lemma}
\label{lmon}
Suppose $M$ is a graded $\mc(n)$-module and $\nabla$ is a flat homogeneous connection. 
Suppose $g:M\r M$ be an isomorphism of $\mc(n)$-modules which is homogeneous
of degree $\ell$. Then 
\beg{eEg}{ \nabla+g^{-1}dg}
is a flat homogeneous connection, and the difference $g^{-1}dg$ has degree $\ell$.
\end{lemma}

\Proof
A direct consequence of Lemma \ref{lmc2}.
\qed

\vspace{3mm}
The connections $\nabla$, \rref{eEg} of Lemma \ref{lmon} will be
said to {\em have the same monodromy}.

\vspace{3mm}

Now we define the degree of a homogeneous flat connection $\nabla$
on $M$. Let us first assume that $M$ is a free graded $\mc(n)$-module.
Let $\phi=(\phi_1,...,\phi_k)$ be a $\mc(n)$-basis of $M$ consisting of
elements of degree $0$. Then we say that the flat homogeneous connection
$$\nabla_\phi:a_1\phi_1+...+a_k\phi_k
\mapsto
(da_1)\phi_1+...+(da_k)\phi_k$$
has degree $0$. An arbitrary flat homogeneous connection $\nabla$ is said to
have degree $0$ when the difference $\nabla-\nabla_\phi$ has degree $0$.
This definition is consistent since by Lemma \ref{lmon},
for two degree $0$ bases $\phi$, $\psi$, $\nabla_\phi-\nabla_\psi$
has degree $0$. 

When $M$ is a projective $\mc(n)$-module, consider a direct summand
embedding $M\subseteq N$ where $N$ is a free $\mc(n)$-module.
We say that a flat connection $\nabla$ on $M$ has degree $0$ when
there exists a commutative diagram
\beg{ehomog+}{\diagram
M\rto^(.3){\nabla}\dto_{\iota}^{\subseteq} & M\otimes_{\mc(n)}\Omega^{1}_{\mc(n)/k}
\dto^{\iota\otimes Id}_{\subseteq}\\
N\rto^(.3){\nabla^\prime} & N\otimes_{\mc(n)}\Omega^{1}_{\mc(n)/k}
\enddiagram
}
where $N$ is a free $\mc(n)$-module and $\nabla^\prime$
is a flat homogeneous connection of degree $0$, and 
$\iota$ is an inclusion of a graded direct summand.

\vspace{3mm}
\begin{lemma}
\label{lhomoga}
This notion of degree $0$ flat homogeneous connection does
not depend on the choice of the graded direct summand inclusion
$\iota$.
\end{lemma}

\Proof
Consider two such direct summand embeddings $\iota:M\r N$,
$\iota^\prime:M\r N^\prime$. Clearly, we may assume
\beg{elhomogai}{N=N^\prime.
}
Let the direct complements of $\iota$, $\iota^\prime$ 
be $K,K^\prime$, respectively. We may assume
\beg{elhomogaii}{K\cong K^\prime
}
(by replacing, if necessary, $N$ with $N\oplus M\oplus K$).
Now assuming \rref{elhomogai}, \rref{elhomogaii},
we can produce a diagram of $\mc(n)$-modules
$$\diagram
M\dto_=\rto^\iota & N\dto_{\cong}^{f}\\
M\rto_{\iota^\prime} & N
\enddiagram
$$
for some graded isomorphism $f:N\r N$ (of degree $0$) of graded $\mc(n)$-modules.
The statement then follows from Lemma \ref{lmon}.
\qed

\vspace{3mm}
\noindent
{\bf Comment:}
We do not know whether every projective $\mc(n)$-module $M$,
or even one endowed with a flat homogeneous connection of
degree $0$, is free. 

\vspace{3mm}

Let $M_i$ be modules over $\Z$-graded commutative rings $R_i$, and suppose we have
homogeneous connections $E_i$ on $M_i$. Then there is a natural 
homogeneous connection $E$ on the $\bigotimes R_i$-module
$\bigotimes M_i$ given by
\beg{emcprod}{E(v_1\otimes ...\otimes v_n)=\cform{\sum}{i=1}{n} v_1\otimes...\otimes v_{i-1}
\otimes E_i(v_i)\otimes v_{i+1}\otimes...\otimes v_n.
}
The same formula also defines a natural homogeneous connection on $M_1\widehat{\otimes} M_2$.

If $R\r S$ is a map of $\Z$-graded commutative rings, $M$ is a graded $R$-module and
we have a homogeneous flat connection $E$ on $M$, then there is a natural ``pushforward"
homogeneous flat connection $E\otimes_R S$ on $M\otimes_R S$, since we have a natural
map 
$$\Omega^{1}_{R/F}\otimes_R S\r \Omega^{1}_{S/F}.$$

\section{Tree functors and tree algebras}

\label{stree}

Let $\Lambda$ be a set (called {\em set of labels}). 
A {\em pre-tree functor} $(\Lambda,M)$ is a system of graded $\mc(n)$-modules
$M_{\lambda_1,...,\lambda_n,\lambda_\infty}$ and graded module isomorphisms
\beg{emcsum}{\diagram\protect
M_{\lambda_{11},...,\lambda_{nm_n},\lambda_{\infty}}\otimes_{\mc(m_1+...+m_n)}
(\mc(m_1)\otimes...\otimes\mc(m_n)\widehat{\otimes}\mc(n))\dto^(.4)\cong \\
\protect\cform{\bigoplus}{\lambda_1,...,\lambda_n}{}
\left(\cform{\bigotimes}{i=1}{n} 
M_{\lambda_{i1},...., \lambda_{im_i},\lambda_i}\right)
\protect
\widehat{\otimes}
M_{\lambda_1,....,\lambda_n,\lambda_\infty}
\enddiagram
}
which must make $M_{\lambda_1,...,\lambda_n,\lambda_\infty}$ a $\Lambda$-sorted $\Z$-graded 
$\mc$-module co-operad. We will also assume that there be a distinguished element $1\in \Lambda$
such that 
\beg{eunit1}{\begin{array}{lll}M_{\lambda_\infty}= & 0 &\text{if $\lambda_\infty\neq 1$}\\
& F=\mc(0) &\text{if $\lambda_\infty=1$}.
\end{array}}
It follows that the co-augmentation \rref{eaug} induce an isomorphism
\beg{eaug1}{\diagram\protect
M_{\lambda_{i+1},...,\lambda_{n},\lambda_\infty}\otimes_{\mc(n-i)}\mc(n)
\rto^\cong & \protect M_{1,...,1,\lambda_{i+1},...,\lambda_{n},\lambda_\infty}
\enddiagram
}
A pre-tree functor is called {\em finite} if all the projective $\mc(n)$-modules
$M_{\lambda_1,...,\lambda_n,\lambda_\infty}$ are finite rank.
A {\em pre-tree algebra} $(\Lambda,M,V,\alpha,\phi)$
consists of a pre-tree functor $M$, a system of $\Z+\alpha_\lambda$-graded
vector spaces $V_\lambda$ such that $\alpha_1=0$, and homomorphisms of $\mc(n)$-modules
\beg{epta}{\phi_{\lambda_1,...,\lambda_n,\lambda_\infty}:M_{\lambda_1,...,\lambda_n,\lambda_\infty}\r
Hom((V_{\lambda_1})\otimes...(\otimes V_{\lambda_n}),
(V_{\lambda_\infty})\otimes_F \mc(n))
}
homogeneous of degree 
\beg{epdeg}{\alpha_\infty-\alpha_1-...-\alpha_n,}
such that if we choose labels $\lambda_{11},...,\lambda_{nm_n},\lambda_\infty$,
and denote, for a graded $\mc(m_1+...+m_n)$-module $X$, put
$$X^\prime :=X\otimes_{\mc(m_1+...m_n)}(\mc(m_1)\otimes...\otimes \mc(m_n)\widehat{\otimes}\mc(n)),
$$
and put
$$\mathcal{M}=\protect\cform{\bigoplus}{\lambda_i}{} 
\left(\cform{\bigotimes}{i=1}{n} 
M_{\lambda_{i1},...\lambda_{in_i},\lambda_{i}}\right)
\protect \widehat{\otimes}M_{\lambda_1,...,\lambda_n,\lambda_\infty},$$
$$\mathcal{V}=\protect\cform{\bigotimes}{\lambda_i}{} \left(
\cform{\bigotimes}{i=1}{n}Hom(\cform{\bigotimes}{j}{} 
V_{\lambda_{ij}},
V_{\lambda_i}\otimes \mc(m_i))\right)\widehat{\otimes}
\protect Hom(\cform{\bigotimes}{i}{}V_{\lambda_i},V_{\lambda_\infty}\otimes \mc(n)),
$$
we have a commutative diagram
\beg{epta1}{\diagram
(M_{\lambda_{11},...,\lambda_{nm_n},\lambda_\infty})^\prime \rto &
\left(Hom(V_{\lambda_{11}}\otimes...\otimes V_{\lambda{nm_n}},\mc(
\cform{\sum}{i=1}{n} m_{i}))
\right)^\prime\\
\mathcal{M}\uto^\cong\rto_{\phi_*\hat{\otimes}\phi_*} &
\mathcal{V},\uto
\enddiagram
}
where the top and bottom row are induced by the appropriate cases of \rref{epta}, the left column
is induced by an appropriate case of \rref{emcsum}, and the right hand column is induced
by composition.

There is also the obvious equivariance axiom and a unitality axiom which asserts that
$1\in M_{\lambda,\lambda}$ maps in
$$Hom(V_\lambda,V_{\lambda}\otimes \mc(1))$$
to $Id\otimes 1$. A pre-tree algebra is called {\em finite} if its pre-tree functor
is finite.

\vspace{3mm}
A {\em tree functor} $(\Lambda,M,E)$
consists of a pre-tree functor
$(\Lambda,M)$, and a system of homogeneous
flat connections 
\beg{efcon}{E_{\lambda_1,...,\lambda_n,\lambda_\infty}}
of degree \rref{epdeg}
on $M_{\lambda_1,...,\lambda_n,\lambda_\infty}$ such that

\begin{enumerate}

\item The connections \rref{efcon} are equivariant with respect to the obvious
action of $\Sigma_n$

\item The connection \rref{efcon} for $\lambda_1=...=\lambda_i=1$ has the same
monodromy as the pushforward of the
connection $E_{\lambda_{i+1},...,\lambda_n,\lambda_\infty}$
along the natural map $\mc(n-i)\r\mc(n)$ induced by co-inserting $\mc(0)$ to the first
$i$ coordinates.

\item \label{efcon3}
The pushforward of the connection 
$$E_{\lambda_{11},...,\lambda_{nm_n}}$$
via the structure map
\beg{efcon31}{\mc(m_1+...+m_n)\r \mc(m_1)\otimes...\mc(m_n)\widehat{\otimes} \mc(n)}
is equal to
$$\cform{\bigoplus}{\lambda_i,...,\lambda_n}{}\left(\cform{\bigotimes}{i=1}{n}
E_{\lambda_{i1},...,\lambda{im_i},\lambda_i}\right)
\widehat{\otimes}
E_{\lambda_1,...,\lambda_n}.
$$

\item The connection $E_{\lambda,\lambda}$ is the pushforward of the connection $dz$
via the map $\mc(1)\r M_{\lambda,\lambda}$.

\end{enumerate}

\vspace{3mm}
A {\em tree algebra} $(\Lambda,M,V,\alpha,\phi,E)$
consists of a pre-tree algebra
$(\Lambda,M,V,\alpha,\phi)$, and a structure $(\Lambda,M,E)$ of a tree functor on its
pre-tree functor $(\Lambda,M)$.

\vspace{3mm}

This data are somewhat redundant. There is no natural choice
of the numbers $\alpha_\lambda$, they are
only determined modulo $1$. Because of that, it is important to define isomorphism
of tree algebras. To this end, there is an obvious notion of
isomorphism of tree functors, and an obvious notion of isomorphism
involving graded isomorphisms of the spaces $V_\lambda$, so all we need to discuss is isomorphism
of tree algebra on the same tree functor and the same
data $V,\phi$. An isomorphism of tree algebras 
$(\Lambda,M,V,\alpha,\phi,E)$, $(\Lambda,M,V,\beta,\phi,F)$
consists of homogeneous isomorphisms
\beg{eiso1}{\diagram
g_{\lambda_1,...,\lambda_n,\lambda_\infty}:M_{\lambda_1,...,\lambda_n,\lambda_\infty}
\rto^(.6)\cong &
M_{\lambda_1,...,\lambda_n,\lambda_\infty}
\enddiagram
}
homogeneous of degree 
$$(\beta_{\lambda_\infty}-\alpha_{\lambda_\infty})
-(\beta_{\lambda_1}-\alpha_{\lambda_1})-...-(\beta_{\lambda_n}-\alpha_{\lambda_n})
$$ 
such that
$$F_{\lambda_1,...,\lambda_n,\lambda_\infty}-E_{\lambda_1,...,\lambda_n,\lambda_\infty}=
g^{-1}_{\lambda_1,...,\lambda_n,\lambda_\infty}dg_{\lambda_1,...,\lambda_n,\lambda_\infty}
$$
and the pushforward of $g_{\lambda_{11},...,\lambda_{nm_n},\lambda_\infty}$ via
\rref{emcsum} is equal to
$$\protect\cform{\bigoplus}{\lambda_1,...,\lambda_n}{}
\bigotimes g_{\lambda_{i1},...., \lambda_{im_i},\lambda_i}\protect\widehat{\otimes} 
g_{\lambda_1,....,\lambda_n,\lambda_\infty}.
$$

\vspace{5mm}

\section{Regularity}

\label{sreg}

\vspace{3mm}

Our treatment of regular connections follows Deligne \cite{d}.
Let $C$ be a smooth algebraic curve over $F$, and let $\overline{C}$
be a smooth projective curve containing $C$. Let $M$ be a (finite-dimensional) algebraic
vector bundle on $C$, and let
\beg{er1}{E: M\r M\otimes_{\mathcal{O}_C} \Omega^{1}_{C/Spec F}
}
be an algebraic connection. We say that $E$ has {\em regular singularities}
if for every smooth projective curve $\overline{C}$ and every embedding $C\subset \overline{C}$,
and every point $x\in \overline{C}-C$, 
$$M\otimes_{\mathcal{O}_C}\mathcal{O}_x$$
has a basis in which the pushforward of the connection $E$ has simple poles. 
A connection $E$ on an $n$-dimensional smooth separated algebraic variety $X$ is said
to have {\em regular singularities} if for every smooth algebraic curve $C$ in $X$,
the restriction of $E$ to $C$ has regular singularities. We will call a tree functor
$(\Lambda,M,E)$ {\em regular} when all the connections $E_{\lambda_1,...,\lambda_n,\lambda_\infty}$
have regular singularities. A tree algebra is called regular when its tree functor
is regular. This notion is clearly invariant under isomorphism of tree algebras.

\vspace{3mm}
Our main goal is to establish a version of the Riemann-Hilbert correspondence for
tree functors and tree algebras over $\C$. To this end, we must define the analytic versions of
our concepts. In effect, we may define $\mc(n)_{an}$ to be the $\Z$-graded ring of all
holomorphic functions $f$ on the configuration space $C(n)$ of $n$ ordered distinct points in
$\C$, homogeneous in the sense that
$$f(\lambda z_1,...,\lambda z_n) =\lambda^k f(z_1,...,z_n);$$
In our grading, the degree of $f$ is $-k$. We would like to
have a $\Z$-graded co-operad structure
\beg{ecoopan0}{
\mc(m_1+...+m_n)_{an}\r (\mc(m_1)_{an}\otimes...\otimes \mc(m_n)_{an})\widehat{\otimes} \mc(n)_{an},}
%and
%\beg{ediffan}{\Omega^{1}_{\mc(n)_{an}/\C}=\mc(n)_{an}\{dz_1,...,dz_n\},
%}
but the difficulty is that expanding an analytic function by a Laurent series
may present elements of arbitrarily low degrees. Replacing $\widehat{\otimes}$
by the product of its bigraded summands of the same total degree, we
do get a $\Z$-graded co-operad, but the resulting objects 
$\mathcal{P}(m_1,...,m_n,n)$
aren't rings,
so we cannot study vector bundles in the sense of finite rank projective modules. 
Our solution is to replace \rref{ecoopan0} by
\beg{ecoopan}{
\mc(m_1+...+m_n)_{an}\r \chi(m_1,...,m_n,n)_{an},}
where the right hand side denotes the ring of all partially defined
holomorphic functions $f$ on 
\beg{ecooprod}{C(m_1)\times...\times C(m_n)\times C(n)}
where, if we denote the coordinates of \rref{ecooprod}
by 
\beg{ecooprod1}{t_{11},...,t_{1m_1},....,t_{n1},...,t_{nm_n}, z_1,...,z_n,}
then for each choice of $z_1,...,z_n$ there exists a locally uniform $\epsilon>0$
such that $f$ is defined on \rref{ecooprod1} when 
\beg{ecooprod2}{||t_{ij}||<\epsilon.
}
Using \rref{ecoopan} and the natural inclusion
\beg{ecooprod3}{\mc(m_1)_{an}\otimes...\otimes \mc(m_n)_{an}\otimes
\mc(n)_{an}\subset \chi(m_1,...,m_n,n)_{an},
} 
we may define analytic (pre)-tree functors and (pre)-tree algebras
in precise analogy with the algebraic definitions. 
We start with $\C^\times$-equivariant holomorphic vector bundles $\Xi$ on $C(n)$.
By this we mean a holomorphic bundle with a holomorphic action of $\C^\times$ on the
total space which is compatible with the $\C^\times$-action on $C(n)$
by
$$\lambda(z_1,...,z_n)=(\lambda z_1,...,\lambda z_n).$$
%For $n=1$, we also additionally assume that the $\C^\times$-action on the fiber over
%$0$ is trivial.
Then the space $M$ of global sections of $\Xi$ is then naturally a graded $\mc(n)_{an}$-module.

\vspace{3mm}
We will called a projective module {\em of finite rank} if it is a direct summand
of a free module on a finite set of generators. 
\vspace{3mm}

\begin{lemma}
\label{lstein}
Let $X$ be a Stein manifold. The global sections functor
defines an equivalence of categories from the category of finite-dimensional
holomorphic vector bundles over $X$ (and holomorphic maps over $Id_X$)
and the category of finite rank projective modules over the ring $Hol(X)$
of holomorphic functions on $X$.
\end{lemma}

\Proof
First we will prove that a holomorphic vector bundle $\xi$ of finite
dimension $n$ on $X$ is always a holomorphic direct summand of 
a finite dimensional trivial vector bundle. First, note that
this is true topologically: $X$ is of the homotopy type of
a $dim(X)$-dimensional CW complex, so by Whitehead's theorem,
the topological
classifying map $\phi$ of $\xi$ factors through
a map $\phi^\prime$ into the $dim(X)$-skeleton
of $BU(n)$:
$$\diagram
X\rto^\phi\drdotted|>\tip^{\phi^\prime}&
BU(n)\\
& BU(n)_{dim(X)}\uto_{\subseteq}
\enddiagram
$$
But $BU(n)_{dim(X)}$ is compact, so the
restriction $\gamma^{\prime}_{n}$ of the
universal $n$-bundle $\gamma_n$ on
$BU(n)$ to $BU(n)_{dim(X)}$ is a direct
summand of a finite-dimensional trivial
vector bundle. Hence, the same is true
for $\xi$ topologically, which we
indicate by the subscript $(?)_{top}$:
\beg{elstein1}{\xi_{top}\oplus\eta_{top}\cong N_{top}.
}
But now the data \rref{elstein1}
may be represented topologically by a 
$GL_n(\C)\times GL_{N-n}(\C)$-principal bundle, so
by Grauert's principle \cite{grauert}, Satz 2,
the data \rref{elstein1} can be represented in
the holomorphic category, which we indicate by
the subscript $(?)_{an}$:
\beg{elstein2}{\xi_{an}\oplus\eta_{an}\cong N_{an}.
}
Applying Grauert's principle again for the groups
$GL_n(\C)$, $GL_N(\C)$, however, we see that
$$\xi_{an}\cong \xi,\; N_{an}\cong N.$$
This shows that the global section functor
in the statement of the Lemma lands in the
category indicated. To show that the functor
is onto on isomorphism classes of objects,
recall that a finite rank projective $Hol(X)$-module
can be constructed from a free module by applying
an idempotent matrix; since a free module always
arises from a trivial bundle, applying the
same matrix on the bundle gives the bundle
corresponding to the finite rank projective
module. 

We now need to prove that the global section
functor is fiathfully full. Since, however, every
object in the source is a holomorphic direct summand of
a trivial finite-dimensional vector bundle,
it suffices to prove that the functor is
faithfully full on the subcategory of finite-dimensional trivial
holomorphic vector bundles, which in turn reduces to
showing that the functor induces bijection of
the set of holomorphic self-maps of the $1$-dimensional
trivial vector bundle to the set of holomorphic
self-maps of the $1$-dimensional free $Hol(X)$-modules.
Obviously, however, both sets are (compatibly) bijective to $Hol(X)$.
\qed

\vspace{3mm}

\begin{corollary}
\label{cstein}
The category of finite-dimensional $\C^\times$-equivariant holomorphic
vector bundles over $C(n)$
and holomorphic $\C^\times$-equivariant
homomorphisms (over the identity on $C(n)$) is equivalent, via the
global sections functor, to the category of finite rank projective
graded $\mc(n)_{an}$-modules and (degree $0$)
homomorphisms
of graded modules.
\end{corollary}

\Proof
In the case $n=1$, the isomorphism class of a $\C^\times$-equivariant bundle
is determined by the representation on the $0$ fiber, which shows that the
functor is a bijection on isomorphism classes of objects. Additionally,
non-equivariantly, the bundles are trivial by Grauert's theorem \cite{grauert},
so in the rank $1$ case,
maps are simply holomorphic functions on $\C$, and therefore homogeneous functions
are simply $z^k$, $k\geq 0$. This is a graded morphism if and only if the $\C^\times$-action
on $0$-fiber of the target is $z^{-k}$ tensored with the action of the $0$-fiber of the
source. This gives the required statement.

In the case $n>1$, the category of $\C^\times$-equivariant holomorphic
bundles and $\C^\times$-equivariant holomorphic
maps is equivalent to the category of holomorphic vector bundles on $C(n)_0$.
Since $C(n)_0$ is a Stein manifold, this is in turn equivalent to the category
of finite rank projective $\mc(n)_0$-modules, which is in turn equivalent to the
category of graded finite type projective $\mc(n)$-modules.

\qed

Also, a holomorphic connection on $\Xi$ gives rise to a connection in the $\mc(n)_{an}$-module
sense
\beg{econnal}{E:M\r M\{dz_1,...,dz_n\};}
note that we have a canonical differentiation
$$\mc(n)_{an}\r \mc(n)_{an}\{dz_1,...,dz_n\}.$$
We may then define flat, homogeneous connections and connections with a given degree
in terms of the algebraic connection \rref{econnal}, and we can mimic the definitions
of the previous section using the space $\chi(m_1,...,m_n,n)_{an}$.
One complication to the compatibility of our notions however
is that $(\mc(m_1)\otimes...\otimes \mc(m_n))\widehat{\otimes}\mc(n)$
is not a subspace of $\chi(m_1,...,m_n,n)_{an}$. Nevertheless, if we
denote by $\chi(m_1,...,m_n,n)$ the intersection of 
$(\mc(m_1)\otimes...\otimes \mc(m_n))\widehat{\otimes}\mc(n)$ and
$\chi(m_1,...,m_n,n)_{an}$ in $P(m_1,...,m_n,n)$,
then the right hand side of \rref{efcon31} can be replaced
by $\chi(m_1,...,m_n,n)$, which allows a comparison.

\vspace{3mm}

We will need to define a regular flat homogeneous connection 
$$E:M^\prime\r M^\prime\{dt_{ij},dz_i\}
$$
on 
a finitely generated
projective
$\chi(m_1,...,m_n,n)$-module $M$ where
$$M^\prime =M\otimes_{\chi(m_1,...,m_n,n)}\chi(m_1,...,m_n,n)_{an}.$$ 
We will say that a flat connection $E$
on $M^\prime$ is regular if each of its solutions has moderate growth.
We say that a solution $g$ has moderate growth if for every 
locally closed smooth algebraic curve $C$ 
in \rref{ecooprod}, every holomorphic embedding $D\subset \overline{C}$
with non-zero derivative at $0\in D$ (where $C$ is a smooth
compactification of $C$) and
every continuous function $\epsilon$ on $D$ such that
for every $x=(t_{ij},z_i)\in D-\{0\}$, $g$ is defined for
$x^\prime=(\delta t_{ij}, z_i)$ whenever $\delta_i\leq \epsilon(x)$,
there exists an $N>0$ such that 
\beg{eregcond}{g(x^\prime)<(\prod\delta_i)^{-N} ||u||^{-N}
}
where $u$ denotes the standard holomorphic coordinate on $D$ and in \rref{eregcond},
$g$ denotes a branch of $g$ on a sector $S$ of $D-\{0\}$ (\cite{d}); we may
define a branch of a solution of $E$ on $S$ as a holomorphic function on
$S$ which satisfies $\nabla g=0$ where $\nabla$ is the pullback
of $E$ to $S$. The multi-valued
section $g$ will also be referred to as {\em regular}.

\vspace{3mm}

Our first version of
Riemann-Hilbert correspondence is the following

\vspace{3mm}

\begin{lemma}
\label{lrh1}
There is an equivalence of categories (canonical up to canonical isomorphism)
between the following categories:
\beg{ecat1}{ \parbox{3.5in}{Finite-dimensional $\mc$-equivariant holomorphic
vector bundles on $C(n)$ with a homogeneous
flat connection of degree $k$}}
\beg{ecat2}{\parbox{3.5in}{Finite rank projective $\Z$-graded $\mc(n)$-modules with a homogeneous regular
flat connection of degree $k$.}}
\end{lemma}

\Proof
When $n=1$, then every then every local system on $C(1)$ is trivial. 
Treating the local system as a free 
$\mc(1)_{an}$-module $M$ with a flat connection $E$,
a priori the free generators of $E$ may not be graded. However, obviously
we have a decreasing filtration $F^k M$ consisting of all elements of
degree $\geq k$. Then, we can consider the associated graded object 
\beg{easgr}{F^k M/F^{k+1}M.}
Looking at the lowest $k$ for which \rref{easgr} is non-trivial,
the fact that $M$ clearly implies that \rref{easgr} is generated by
homogeneous elements of degree $k$. Further, the connection must be
$0$ on these generators by homogeneity. Consider the free graded submodule
generated by these elements by $M_0$. Then $M/M_0$ is also a free module
(the category of finite dimensional free $\mc(1)_{an}$-modules with
flat connection is equivalent to the category of finite dimensional
$\C$-vector spaces), so we may repeat the argument with $M$ replaced
by $M/M_0$ to show by induction that $M$ is free as a graded $\mc(1)_{an}$-module
with generators annihilated by the connection. Since clearly an analogous
argument applies to the algebraic category, the statement follows.

Assume now $n>1$. Then any $z_i$ defines an isomorphism of vector
spaces from $\mc(n)_k$ to $\mc(n)_{k-1}$.
The category of graded $\mc(n)$-modules and
morphisms of degree $0$ is therefore
equivalent to the category of $\mc(n)_0$-modules.
The statement for flat connections is that flat connections o
Similar statements are also valid for the
corresponding analytic categories.
But now $\mc(n)_0$ is in fact the coefficient ring of
a smooth affine variety $C(n)_0$, which comes
with an embedding into $\P^{n-1}$ with homogeneous
coordinates $z_1,...,z_n$ (in fact, the embedding 
factors through the copy of $\A^{n-1}$ which is 
the complement of the locus of $z_1-z_2$). 
The statement for connections is that flat connections (algebraic
or analytic) on $C(n)_0$ correspond precisely to
homogeneous connections on $\mc(n)$ resp. $\mc(n)_{an}$
of degree $0$: This is because 
$$\Omega^{1}_{\mc(n)_0/F}=\{a_1 \frac{dz_1}{z_1}+...+a_n\frac{dz_n}{z_n}\in (\Omega^{1}_{\mc(n)/F})_0
| a_1+...+a_n=0\}.$$
Therefore, for homogeneous flat connections of degree $0$, the
result follows from Theorem 5.9 of \cite{d}, applied to
the variety $C(n)_0=Spec(\mc(n)_0)$. The case of general degree $k$ can be reduced to the
case of degree $0$ by subtracting an algebraic connection of degree $k$.
\qed

%\vspace{3mm}
%The same argument, in fact, implies the following generalization, which we
%will also need:

%\begin{lemma}
%\label{lrh2}
%Let $n_1,...,n_p$ be natural numbers.
%Then there is an equivalence of categories (canonical up to canonical isomorphism)
%between the following categories:
%\beg{ecat3}{ \parbox{3.5in}{Finite rank projective $\Z$-graded $(\mc(n_1)\otimes...
%\otimes \mc(n_p))_{an}$-modules with a homogeneous
%flat connection of degree $k$}}
%\beg{ecat4}{\parbox{3.5in}{Finite rank projective $\Z$-graded $\mc(n_1)\otimes...
%\otimes \mc(n_p)$-modules with a homogeneous regular
%flat connection of degree $k$.}}
%(The degree of a homogeneous connection on $\mc(n_1)\otimes...\otimes \mc(n_p)$
%is defined as the degree of its image in $\mc(n_1+...+n_p)$ via the natural localization inclusion.)
%\end{lemma}

%\qed

\vspace{3mm}
Our main interest is in the following statement:

\begin{theorem}
\label{trh}
There is an equivalence of categories (canonical up to canonical isomorphism)
between the category of finite analytical tree functors (resp. algebras) and finite regular
tree functors (resp. algebras).
\end{theorem}

\Proof
Given Lemma \ref{lrh1}, and the fact that the two connections whose
isomorphism we are seeking are obviously regular, 
the statement amounts to asserting that an isomorphism 
in
\beg{ehatan}{\chi(m_1,...,m_n,n)_{an}}
between
two graded regular 
connections on
\beg{ehat1}{\chi(m_1,...,m_n,n)
}
of the same degree is algebraic. However, 
this amounts to saying that every regular function
in \rref{ehatan} is in \rref{ehat1}.
This can be shown as follows: consider a regular function in
\rref{ehatan}. Then considering $g$ as a function of the $t_{ij}$'s
for fixed $z_i$, expand in the total degree of all the $t_{ij}$'s
for a fixed $i$.
Then the coefficients of fixed total degree $d_i$ in the $t_{ij}$'s
for each $i$ are obviously elements of 
$$\mc(m_1)_{d_1}\otimes...\otimes\mc(m_n)_{d_n}$$
and further each of the numbers $d_i$ where the coefficient is non-zero
is bounded below by a bound $N(z_1,...,z_n)$.
Since this bound, however, is locally constant, it must be in effect constant because
the function involved are analytic. 

Now we claim that for each given assortment of degrees $d_1,...,d_n$,
the corresponding component 
\beg{egcomp1}{g_{d_1,...,d_n}}
of $g$ is in effect an element of
$$\mc(m_1)\otimes...\otimes\mc(m_n)\otimes\mc(n).$$
In fact, select the lowest $(d_1,...,d_n)$ (say, lexicographically) such that
\rref{egcomp1} is non-trivial. Then by \rref{eregcond}, all summands of higher
degree can be neglected, and \rref{egcomp1} must be regular. Subtracting this component,
we may show by induction that all of the components are regular. 

\qed

\vspace{3mm}

\begin{definition}
\label{dvt}
Recall the definition of vertex tensor category
of Huang and Lepowsky \cite{hl}, Definition 4.1.
We call a vertex tensor category
{\em semisimple} if there are finitely many
objects (called {\em irreducible objects})
whose endomorphism groups are $\mc$ and such that
there are no nonzero morphisms between non-isomorphic irreducible
objects, and
every object is isomorphic to a direct sum of irreducible objects,
and the unit object is irreducible.
\end{definition}

\vspace{3mm}

\begin{theorem}
\label{thuang}
There is a `realization' functor (canonical up to natural equivalence)
from the category of vertex tensor categories and isomorphisms, 
to the category of analytic tree algebras and isomorphisms, and consequently, by Theorem
\ref{trh}, to the category of tree algebras and isomorphisms. 
\end{theorem}

\vspace{3mm}

\noindent
{\bf Comment:} By the work of Lepowsky and Huang, vertex algebras 
which satisfy certain `rationality' conditions supply examples of vertex
tensor categories in the sense of \cite{hl}. A good survey is Huang \cite{hrep}.
It should also be noted that the present result overlaps with Huang's
construction \cite{h7} of genus $0$ correlation functions for modules of
vertex algebras. The major point of interest of our result is that it
extends to the
algebraic category of tree algebras over $\C$. It should be noted that we model only a part of the
structure of vertex tensor
category in our axioms (e.g. we do not treat the
conformal element = energy-momentum tensor), which is one of the
reasons why we do not get 
an equivalence of categories in our statement.

\vspace{3mm}

\noindent
{\bf Proof of Theorem \ref{thuang}:} We will study the definition
of vertex tensor category \cite{hl}, Definition 4.1.
First of all, because we
do not treat conformal element data, we restrict 
attention to
the subspaces $K(n)_0$ of the moduli spaces $K(n)$ (see Huang \cite{huangbook}, p.65)
where the tube functions are the identity, and the 
scaling constant is $1$.
In this setting, there is a canonical splitting of the determinant bundle,
so we get a canonical map 
\beg{ekn0}{\psi:K(n)_0\r \widetilde{K}^c(n).}
Next, for a semisimple vertex tensor
category $\mathcal{V}$, the set of labels $\Lambda$ is the 
set of representatives of isomorphism classes of irreducible
objects (we shall also write $V_\lambda=\lambda$), and $1$ corresponds to 
the 
unit object \cite{hl}, Definition 4.1. (3). 
Now for 
and $Q\in K(2)_0$, and irreducible objects $V_\lambda$, 
$V_\mu$, we consider their tensor product \cite{hl}, Definition 4.1 (1)
\beg{ehltensor}{W:=V_\lambda \boxtimes_{\psi(Q)}V_\mu.
}
Then by our definition of semisimple vertex tensor category, we have
a unique decomposition up to isomorphism
\beg{eqd1}{W=\cform{\bigoplus}{\nu\in \Lambda}{}W_\nu
}
where, non-canonically,
\beg{egd2}{W_\nu\cong N\otimes V_\nu}
for a finite-dimensional complex vector space $N$.
Further, by our assumptions, non-canonically, we have
\beg{egd3}{Aut(W_\nu)\cong GL(N).
}
Property (6), together with axioms (1) and (2)
of Definition 4.1 of \cite{hl} imply that
\rref{egd3} define a principal smooth $GL(N)$-bundle with 
flat connection on the space $C(2)$. We let $M_{\lambda,\mu,\nu}$
be the dual of the
associated vector bundle (which the flat connection automatically
makes analytic). The correlation function $\phi_{\lambda,\mu,\nu}$
(the analytic version of \rref{epta})
then follows from the universal intertwining operator
from $V_\lambda\otimes V_\mu$ to \rref{ehltensor}.

General correlations functions with an arbitrary number
of arguments are then produced in
an analogous way by iterating
the tensor product \rref{ehltensor}. Co-operad
associativity resp. unitality resp. equivariance follow from property 
(4) resp. (7) resp. (5) of \cite{hl}, Definition 4.1.
Vertex algebra unitality follows from property (8).
\qed

\vspace{3mm}
\noindent
{\bf Example:} The chiral WZW model. Let $\mathfrak{g}$ be a
finite-dimensional simple Lie algebra over $\C$.
In \cite{hl1}, Huang and Lepowsky construct a vertex
tensor category of finite sums of irreducible lowest 
weight-modules $L(k,\lambda)$ over the quotient $L(k,0)$
of the affine vertex algebra $M(k,0)$ by its maximal
ideal (=maximal proper graded submodule) for $k=0,1,2,...$.
By Theorem \ref{thuang}, this gives rise to a tree algebra $T$
over $\C$. In effect, we have the following refinement:

\vspace{3mm}
\begin{theorem}
\label{thuang1}
Let $\mathfrak{g}$ be defined over $\C\supseteq k\supseteq \Q$.
Then the tree algebra corresponding to the Huang-Lepowsky construction
can be defined over $k$.
\end{theorem}

\Proof
The key point is to study the so called Knizhnik-Zamolodchikov
equations \cite{hl1}, (2.14),
which define the desired homogeneous flat connection on the
tree functor. The connection is
defined on the trivial $\mc(n)$-module
\beg{ekz*}{N_{\lambda_1,...,\lambda_n,\lambda_\infty}=
L(\lambda_1)\otimes...\otimes L(\lambda_n)
}
where $L(\lambda)$ is the summand of lowest degree
of $L(k,\lambda)$. The flat connection defined
by the KZ-equations maps 
$$f:L(\lambda_1)\otimes...\otimes L(\lambda_n)$$
to
\beg{ekz+}{df-\frac{1}{k+h^{\vee}}\cform{\sum}{p\neq\ell}{}
\frac{1}{z_\ell-z_p}\cform{\sum}{i}{}
f(Id\otimes...\otimes g^i\otimes...\otimes g_i\otimes...\otimes Id)dz_\phi}
where on the right hand side, $(g^i)$ and $(g_i)$ are dual bases of
$\mathfrak{g}$ with respect to the Killing form, and are
inserted at the $\ell$'th and $p$'th coordinate, respectively.

Manifestly, the connection \rref{ekz+} is defined over $\Q$. 
The tree functor is actually a direct summand of \rref{ekz*}.
It corresponds to $\mathfrak{g}$-equivariant maps from the
lowest weight summand of 
$$(L(k,\lambda_1)\boxtimes...\boxtimes L(k,\lambda_n))_Q$$
to $L(k,\lambda_\infty)$ where $Q$ encodes the moduli
data \cite{hl}, which, in our case, is just the $n$-tuple
of points $(z_1,...,z_n)$. By definition, all this 
is defined over $k$.
\qed

\end{document}